\documentclass[11pt]{amsart}
\usepackage[T1]{fontenc}
\usepackage[ansinew]{inputenc}
\usepackage{amssymb}
\newtheorem{thm}{Theorem}[section]

\newtheorem{corollary}[thm]{Corollary}
\newtheorem{lem}[thm]{Lemma}
\newtheorem{probl}{Problem}

\newtheorem{proposition}[thm]{Proposition}
\theoremstyle{definition}
\newtheorem{definition}[thm]{Definition}
\theoremstyle{remark}
\newtheorem{remark}{Remark}[section]

\theoremstyle{procedure}

\theoremstyle{Example}

\theoremstyle{Assumption}

\newtheorem{asmp}[thm]{Assumption}

\newcommand{\ad}{\mbox{ad}}

\newcommand{\tu}{\tilde{u}(\cdot)}
\newcommand{\txe}{\tilde{x}}
\newcommand{\tue}{\tilde{u}}

\newcommand{\tx}{\tilde{x}(\cdot)}
\newcommand{\tps}{\tilde{\psi}}
\newcommand{\spn}{\mbox{Span }}

\theoremstyle{remark}
\newtheorem{prethm}{Theorem}
\theoremstyle{remark}
\newtheorem{predef}[prethm]{Definition}
\theoremstyle{remark}
\newtheorem{precol}[prethm]{Corollary}
\theoremstyle{remark}
\newtheorem{preprop}[prethm]{Proposition}
\theoremstyle{remark}
\newtheorem{prerem}[prethm]{Remark}

\begin{document}
\title{Controlling Multiparticle System on a Line. I}%
\author{Andrey
V. Sarychev$^1$}%
\address{$^1$ DiMaD,
University of Florence, Italy}%
\email{asarychev@unifi.it}%
\date{}
\begin{abstract}
We  study a classical multiparticle system (such as Toda lattice)
whose dynamics we intend to control by forces applied to few
particles of the system. Various problem settings, typical for
control theory are posed for this model; among those: studying
accessibility and controllability properties, structure properties
and feedback linearization of respective control system,
time-optimal relocation of particles. We obtain complete or partial
answers to the posed questions; criteria and methods of geometric
control theory are employed. In the present part I we consider
nonperiodic multiparticle system. In the forthcoming Part II we
address controllability issue for multiparticle system subject to
periodic boundary conditions. That study requires an extension and
refinement of known methods of geometric control.
\end{abstract}
\maketitle

\vspace*{2mm}

{\small\bf AMS Subject Classification: } %

\section{Introduction}
\label{intro} \markboth{A.V.Sarychev}{Controlling Multiparticle
System}

Consider classical system of $n$ interacting particles
$\mathcal{P}_1, \ldots , \mathcal{P}_n$ moving on a  line with only
neighboring particles involved in interaction. Let $q_k$ be the
coordinate of the $k$-th particle and $p_k$ - its impulse.

We assume the potential of the interaction to be
$$\Phi(q_1-q_{2})+\Phi(q_2-q_{3})+\cdots +\Phi(q_{n-1}-q_n),$$
where $\Phi: \mathbb{R} \rightarrow \mathbb{R}$ is real analytic
 \footnote{In fact most part of the results below would be valid for $C^\infty$-smooth $\Phi$,
 but the reasoning in the real analytic case is less technically
 involved.},
 bounded
below function
\begin{equation}\label{grof}
\lim_{y \rightarrow +\infty}\Phi(y)=+\infty.
\end{equation}

The dynamics of such system of particles is described by the
Hamiltonian system of equations with the Hamiltonian
$$H(q,p)=\frac{1}{2}\sum_{k=1}^np_k^2+\sum_{j=1}^{n-1}\Phi(q_j-q_{j+1}).$$
In coordinates $q_k,p_k$ the equations of system are
\begin{eqnarray}
\dot{q}_k=\frac{\partial H}{\partial p_k}=p_k, \ k=1, \ldots ,n, \label{dqq}\\
\dot{p}_k=-\frac{\partial H}{\partial
q_k}=\phi(q_{k-1}-q_k)-\phi(q_{k}-q_{k+1}), \ k=2, \ldots , n-1, \label{dpq}\\
\dot{p}_1=-\phi(q_{1}-q_{2}), \ \dot{p}_n=\phi(q_{n-1}-q_{n})
\label{dp1n}
\end{eqnarray}
where $\phi=\Phi '$ is the derivative of $\Phi$. It is natural to
assume $$\lim_{y \rightarrow -\infty}\phi(y)=0,$$
 the interaction
decreases to zero, when the distance between particles tends to
infinity. Under this additional assumption we can adapt the
equations (\ref{dp1n}) to the form (\ref{dpq}), introducing
fictitious particles $\mathcal{P}_0$ and $\mathcal{P}_n$ on which we
impose  boundary conditions
$$q_0=-\infty, \ q_{n+1}=+\infty.$$

Our main goal will be  controlling the location and the momenta of
the particles in situation where the control tools are limited, that
is control are forces applied only to few particles of the system.

We will study two cases: {\it single forced  multiparticle system}
with a controlled force acting only on the particle $\mathcal{P}_1$
(or on $\mathcal{P}_n$), {\it double forced multiparticle system}
with controlled force   applied to the particles $\mathcal{P}_1$ and
$\mathcal{P}_n$.

It turns out that controlled multiparticle system provides a model
example for application of the methods of geometric control theory.
Below we demonstrate their effectiveness  for solution of various
control-theoretic problems for this model. In Section~\ref{liacc}
start studying  the Lie structure of the multiparticle system,
verifying full-dimensionality of its orbits and zero-time orbits and
establishing its strong accessibility, whenever the system is
controlled by single force applied to either $\mathcal{P}_1$ or
$\mathcal{P}_n$. We establish property of global controllability for
double forced system in Section~\ref{dic}. The subsequent study of
the Lie structure of single and of double forced multiparticle
systems Section~\ref{feelin} show that in many aspects these systems
behave like linear ones. This is validated by result on their local
feedback linearizability. The linear-like structure reveals again
when we study in Section~\ref{toc} time-optimal particle relocation
problem by means of {\it constrained} controls. We prove that the
corresponding time-optimal controls are bang-bang, i.e. admit their
values at extreme points of the rectangle which constrains the
control parameters.
 The number of switchings between these extreme
points is proved to be uniformly bounded for trajectories evolving
on a fixed compact of state space.

Another model, which  we shall study  in the forthcoming Part II of
the publication,  is multiparticle system under periodic boundary
conditions
$$q_0=q_n, \ q_{n+1}=q_1.$$
 It is known that the dynamics
of the periodic and nonperiodic Toda lattices are completely
different. Also the Lie structures of nonperiodic and periodic
controlled multiparticle systems differ substantially: the latter is
 far from being linear-like, what we establish for the nonperiodic system below. In Part II are going to study
controllability properties of controlled periodic multiparticle
system.

\section{Single-forced multiparticle system; Lie structure and accessibility property}
\label{liacc}

We introduce  control $u_1(t)$ which is time-varying  force applied
to the particle $\mathcal{P}_1$ of the multiparticle system. We
obtain then for the (momentum) variable $p_1$ the equation
\begin{equation}\label{ceqb1}
\dot{p}_1=-\phi(q_1-q_2)+u(t).
\end{equation}
The equations  (\ref{dpq}) and the second equation in (\ref{dp1n})
remain unchanged.

The  controlled multiparticle system is a particular case
single-input control-affine system of the form
\begin{equation}\label{caf}
    \dot{x}=f(x)+g^u(x)u, \ u \in \mathbb{R},
\end{equation}
where $x=(q,p)=(q_1, \ldots , q_{n},p_1, \ldots , p_{n}) \in
\mathbb{R}^{2n}$, the uncontrolled vector field $f$ - the {\it
drift} -  and the {\it controlled vector field} $g$ are defined as
\begin{eqnarray}\label{vef}
g^u=\frac{\partial}{\partial p_1}, \
f=\sum_{k=1}^{n}p_k\frac{\partial}{\partial q_k}-
\phi(q_1-q_2)\frac{\partial}{\partial p_1}+ \\
+\sum_{k=2}^{n-1}
\left(\phi(q_{k-1}-q_k)-\phi(q_{k}-q_{k+1})\right)\frac{\partial}{\partial
p_k}+\phi(q_{n-1}-q_n)\frac{\partial}{\partial p_n}. \nonumber
\end{eqnarray}

In this Section we start studying  Lie structure of this
single-input control-affine system by establishing {\it
accessibility property} - full-dimensionality of its orbits and
attainable sets. Exact definitions and needed criteria are provided
in the following Subsection.

\subsection{Preliminaries}

 \subsubsection{Vector fields, Lie brackets}
Real analytic vector field in $\mathbb{R}^N$ is an analytic map $x
\mapsto F(x) \in T_x\mathbb{R}^N \simeq \mathbb{R}^N$.

Any vector field $F$ defines derivation $\hat{F}$ of the algebra of
analytic functions on $\mathbb{R}^N$ and vice versa. The commutator
of two derivations $\hat{F}^1,\hat{F}^2$ is again a derivation, and
the  corresponding vector field is called the {\it Lie bracket}
$[F^1,F^2]$ of $F^1,F^2$. The operation $[\cdot , \cdot]$ defines
structure of Lie algebra in the space of vector fields. In
coordinates it is calculated as
$$[F^1,F^2]=DF^2 F^1 - DF^1 F^2,$$
where $DF$ stays for the Jacobian matrix of $F$.

For a vector field $F$ we consider the operator $\ad F$, which acts
in the Lie algebra of vector fields in $\mathbb{R}^{2n}$: $\ad F
F^1=[F,F^1] $. The iterations of this operator are denoted $\ad^j
F$: $\ad^j F F^1=[F, \ad^{j-1}F F^1]$.  This operator is a
derivation of the Lie algebra; it satisfies the Leibniz rule:
  $$\ad f [g,h]=[\ad f g, h]+[g, \ad f h]$$
which is equivalent to the Jacobi identity of the Lie algebra.

\subsubsection{Lie envelope, zero-time ideal}
\label{lienv} Below we introduce the needed notions and formulate
results for the class of control-affine systems; readers may consult
 \cite{Ju,ASkv} for the same material in more general context.

Consider a control-affine system
\begin{equation}\label{afgen}
\dot{x}=f(x)+G(x)u=f(x)+\sum_{i=1}^rg^i(x)u_i, \  u=(u_1, \ldots ,
u_r) \in U,
\end{equation}
 where $f, g^1, \ldots , g^r$ are real analytic vector fields on $\mathbb{R}^N$.
 We assume  the set $U$ of control parameters to contain the
 origin $0_{\mathbb{R}^r}$ in its interior.

Let $\mbox{Lie}\{f,G\}$ be the Lie algebra generated by $f,g^1,
\ldots , g^m$, and $\mathcal{I}^0\{f,G\}$ be its Lie ideal generated
by $g^1, \ldots , g^m$. We will call them the {\it Lie envelope} and
{\it zero-time ideal} of the control system respectively.

\subsubsection{Orbits. Orbit Theorem}
\label{orbeth}

Substituting constant controls $u^j=(u^j_1, \ldots , u^j_r)$ into
the right-hand side of \eqref{afgen} we obtain vector fields
$f^{u^j}$ which generate corresponding flows $e^{t_jf^{u^j}}$.
Acting by the compositions
\begin{equation}\label{calP}
P=  e^{t_1f^{u^{j_1}}} \circ \cdots \circ e^{t_Nf^{u^{j_N}}}, \ t_1,
\ldots , t_N \in \mathbb{R},
\end{equation}
onto a given point $\tilde{x}$ we get  an {\em orbit}
$\mathcal{O}_{\tilde{x}}$ of the control system \eqref{afgen} from
$\tilde{x}$. Requiring in addition $t_1+ \cdots +t_N=0$ at
 the right-hand side of \eqref{calP}, we get {\it zero-time
orbit $\mathcal{O}^0_{\tilde{x}}$ of the system}.

Acting by a diffeomorphism $P^0=e^{t_0f^{u^0}}$ (or by a composition
\eqref{calP} of such diffeomorphisms) onto zero-time orbit
$\mathcal{O}^0_{\tilde{x}}$ we obtain zero-time orbit
$\mathcal{O}^0_{y}$, where $y=P^0(x)$.

The orbits and  zero-time orbits possess regular structure.

\begin{prethm}[Orbit Theorem; Nagano-Stefan-Sussmann, \cite{ASkv}]
\label{nss}  An orbit $\mathcal{O}_{\tilde{x}}$ and  zero-time orbit
$\mathcal{O}^0_{\tilde{x}}$ of the control system \eqref{afgen} are
immersed submanifolds of $\mathbb{R}^N$. The tangent space to the
orbit  $\mathcal{O}_{\tilde{x}}$ at a point $x \in
\mathcal{O}_{\tilde{x}}$ coincides with the evaluation
$\mbox{Lie}_{x}\{f,G\}$ of the vector fields from
$\mbox{Lie}\{f,G\}$ at $x$; the tangent space to the zero-time orbit
$\mathcal{O}^0_{\tilde{x}}$ at a point $x \in
\mathcal{O}^0_{\tilde{x}}$ coincides with the evaluation
$\mathcal{I}^0_{x}$ of the vector fields from $\mathcal{I}^0$ at $x.
\ \Box$
\end{prethm}

\begin{precol}\label{crank}
The dimensions $d(x)=\dim \mbox{Lie}_{x}, \ d^0(x)=\dim
\mathcal{I}^0_{x}$ are  constant  along any  orbit and  zero-time
orbit of the system respectively. Obviously $d^0(x) \leq d(x) \leq
d^0(x)+1. \ \square$
\end{precol}

\subsubsection{Attainable sets. Accessibility
property} \label{acprop}

Involving only those compositions \eqref{calP}, where $t_j$ are
nonnegative, and acting by them on a given point $\tilde{x}$ we
obtain {\it positive orbit} or {\it attainable set}
$\mathcal{A}_{\tilde{x}}$ of the system \eqref{afgen} from
$\tilde{x}$. If we pick $T>0$ and require in addition  $t_1 + \cdots
+t_N=T$, or respectively, $t_1 + \cdots +t_N \leq T$,  then we
obtain time-$T$ (respectively time-$\leq T$) attainable set
$\mathcal{A}^T_{\tilde{x}}$ (respectively $\mathcal{A}^{\leq
T}_{\tilde{x}}$).

Obviously $\forall \tilde{x} \ \forall T>0$:
$$ \mathcal{A}^T_{\tilde{x}} \subset \mathcal{A}^{\leq T}_{\tilde{x}} \subset \mathcal{A}_{\tilde{x}}
\subset \mathcal{O}_{\tilde{x}}.$$ Besides
$$\mathcal{A}^T_{\tilde{x}} \subset e^{Tf}\left( \mathcal{O}^0_{\tilde{x}}\right).$$

It turns out that $\mathcal{A}_{\tilde{x}}$ and
$\mathcal{A}^T_{\tilde{x}}$ are 'massive' subsets of
 $\mathcal{O}_{\tilde{x}}$ and of $e^{Tf}\left(
 \mathcal{O}^0_{\tilde{x}}\right)$ respectively.

\begin{prethm}[Krener]
\label{krener}  Attainable set $\mathcal{A}^{\leq T}_{\tilde{x}}$
and $\mathcal{A}^T_{\tilde{x}}$ possess nonvoid relative interiors
in  $\mathcal{O}_{\tilde{x}}$ and  in  $e^{Tf}\left(
\mathcal{O}^0_{\tilde{x}}\right)$) respectively. The sets
$\mathcal{A}_{\tilde{x}},\mathcal{A}^T_{\tilde{x}}$ are contained in
the closures of their relative interiors. $\Box$
\end{prethm}

By virtue of Theorems~\ref{nss} and \ref{krener} the sets
$\mathcal{A}^{\leq T}_{\tilde{x}}$ (respectively
$\mathcal{A}^T_{\tilde{x}}$) posses absolute interior whenever $\dim
\mbox{Lie}_{x}=N$ (respectively $\dim \mathcal{I}^0_{x}=N$).

In these cases the control system \eqref{afgen} is said to possess
{\it accessibility property} (respectively  {\it time-$T$
accessibility property}) from a point $x$.

We will be interested in stronger property of  {\it global
controllability}.
\begin{predef}
The system is globally controllable if
$\mathcal{A}_{\tilde{x}}=\mathbb{R}^N. \ \Box$
\end{predef}

It is immediate to see that global controllability implies
accessibility property.

\begin{preprop}
$$\mathcal{A}_{\tilde{x}}=\mathbb{R}^N
\Longrightarrow \mathcal{O}_{\tilde{x}}=\mathbb{R}^N
\Longleftrightarrow \dim \mbox{Lie}_{x}=N \Longrightarrow
\mbox{accessibility}. \ \Box$$
\end{preprop}

The inverse implication is not valid; $\mathcal{A}_{\tilde{x}}$ is
'often' a proper subset of $\mathcal{O}_{\tilde{x}}$. A sufficient
criterion for coincidence of these two sets is discussed in
Subsection~\ref{prerecur}.

\subsection{Lie envelope and accessibility property for single-forced multiparticle system}

Coming back to the control system \eqref{caf}-\eqref{vef} we are
going to calculate dimension of its  orbits  and establish
accessibility property. The results we obtain remain valid whenever
control is applied to $\mathcal{P}_n$ instead of $\mathcal{P}_1$.

By virtue of the  criteria  formulated in the
Subsubsections~\ref{acprop} and \ref{orbeth} both properties can be
derived from the following technical proposition.

\begin{proposition}
\label{lier}  The dimension $\dim \mathcal{I}^0_{x}$ equals $2n$ at
each point $x \in \mathbb{R}^{2n}. \  \Box$
\end{proposition}

This Proposition would follow immediately from the following Lemma,
which provides more information on the Lie structure of
\eqref{caf}-\eqref{vef}.

\begin{lem}
\label{complam} For each $k \geq 0$ the distributions
\begin{equation}\label{xlam}
x \mapsto \Lambda^m_x=\spn \{(\ad^k f g^u)(x), k=0, \ldots , m-1\},
\ \Lambda^0=\{0\},
\end{equation}
meet the relations
\begin{equation}\label{lamk}
\Lambda^{2k}_x \subseteq  \spn\left\{\frac{\partial}{\partial
p_s},\frac{\partial}{\partial q_s}\left|\  1 \leq s \leq
k\right.\right\},\ \Lambda^{2k+1}_x \subseteq
\Lambda^{2k}+\spn\left\{\frac{\partial}{\partial p_{k+1}}\right\}
\end{equation}
with equalities in \eqref{lamk} holding at a generic point of a
zero-time orbit of the system (\ref{caf})-\eqref{vef}. $\Box$
\end{lem}

Assuming validity of the conclusion of the Lemma we pick any point
$x \in \mathbb{R}^{2n}$ and consider the corresponding zero-time
orbit $\mathcal{O}^0_x$. At a generic point of this orbit
$\mathcal{I}^0_{x} \supset \Lambda^{2n}$ and hence $\dim
\mathcal{I}^0_{x}=2n$-dimensional. Then its (constant) dimension is
$2n$ at each point of $\mathcal{O}^0_x$.

An immediate corollary of the Proposition \ref{lier}  is the {\it
accessibility property}.

\begin{thm}
The multi-particle system,  controlled by a single control (force),
applied to either $\mathcal{P}_1$ or $\mathcal{P}_n$, possesses for
any $T>0$, time-$T$ accessibility property, and the set
$\mathcal{A}^T_{\tilde{x}}$ of positions $q$ and momenta $p$ of the
particles attainable from $\tilde{x}=(\tilde{p},\tilde{q})$ in any
time $T>0$ has an interior, which is dense in
$\mathcal{A}^T_{\tilde{x}}. \ \Box$
\end{thm}

A relevant question is what happens with accessibility when the
controlled force is applied to an 'intermediate' particle
$\mathcal{P}_j, \ j \neq 1,n$. In this case the Lie structure is not
so regular as the one defined by \eqref{xlam}-(\ref{lamk}). In fact
for a generic $\phi$ the Lie rank is complete and the system
possesses the accessibility property. Still for special choice of
$\phi$ the system may possess low-dimensional orbits and therefore
lack the accessibility property. We provide corresponding example in
the forthcoming Part II of the publication where we study this and
other issues for multiparticle system periodic boundary conditions.

{\it Proof of Lemma~\ref{complam}}. First, note that
\begin{eqnarray}
\left[f,\frac{\partial}{\partial
p_s}\right]=-\frac{\partial}{\partial q_s}, \label{bfq}\\
 \left[f,\frac{\partial}{\partial
q_s}\right]=-\phi'(q_{s-1}-q_s)\left(\frac{\partial}{\partial
p_{s-1}}-\frac{\partial}{\partial p_s}\right) +\\
+\phi'(q_s-q_{s+1})\left(\frac{\partial}{\partial p_s}-
\frac{\partial}{\partial p_{s+1}}\right), \nonumber
 \ s=2, \ldots ,n-1,
\label{bfp}\\
\left[f,\frac{\partial}{\partial
q_1}\right]=\phi'(q_1-q_2)\left(\frac{\partial}{\partial
p_1}-\frac{\partial}{\partial p_2}\right), \label{fq1} \\
\left[f,\frac{\partial}{\partial q_n}\right]=
-\phi'(q_{n-1}-q_n)\left(\frac{\partial}{\partial
p_{n-1}}-\frac{\partial}{\partial p_n}\right). \label{fqn}
\end{eqnarray}

Now we proceed by induction on $k$ proving (\ref{lamk}) and
verifying at the same time, that
\begin{eqnarray}\label{2k1}
\ad^{2k-1}fg=(-1)^{k-1}\mu_{k-1}(q)\frac{\partial}{\partial p_{k}} \
(\mbox{mod
}\Lambda^{2k-2}),\\
\label{2k2} \ad^{2k}fg= (-1)^{k-1}\mu_{k-1}(q)
\frac{\partial}{\partial q_{k}} \ (\mbox{mod }\Lambda^{2k-1}).
\end{eqnarray}
where $ \mu_k(q)=\prod_{j=1}^k \phi'(q_j-q_{j+1})$ and we assume
$\mu_k=1$ for $k=0$.

 For $\Lambda^1,\Lambda^2$ formulae (\ref{delk}) are valid, while formulae (\ref{2k1})-(\ref{2k2})
 are trivial.

Let $\Lambda^{2k}$ be the distribution defined by (\ref{xlam}) with
$m=2k$.  Our induction assumption is that \eqref{2k1} and
(\ref{lamk}) are valid for $\Lambda^{2k}$.  According to
\eqref{lamk} the vector fields $\ad^\ell fg^u$ with $\ell < 2k$ can
be represented as $\sum_{s=1}^k\alpha_s(x)\frac{\partial}{\partial
q_s}+ \beta_s(x)\frac{\partial}{\partial p_s}$
 To evaluate
$[f,\Lambda^{2k}]$ we consider the Lie bracket $\left[f,
\sum_{s=1}^k\alpha_s(x)\frac{\partial}{\partial q_s}+
\beta_s(x)\frac{\partial}{\partial p_s}\right]$ and conclude by
(\ref{bfq}),(\ref{bfp}) that its values are contained in
$$\spn\left\{\frac{\partial}{\partial p_j},\frac{\partial}{\partial
q_s}\left|\ j=1, \ldots , k+1;\ 1 \leq s \leq k\right.\right\}.$$

On the other side by induction hypothesis
$$\ad^{2k-1}fg= \sum_{s=1}^k\alpha_s(x)\frac{\partial}{\partial q_s}+
\beta_s(x)\frac{\partial}{\partial p_s},$$ with
$\alpha_k=(-1)^k\mu_k(q)$, being nonvanishing at a generic point.

The following equalities hold modulo $\Lambda^{2k}$ at chosen point
of the orbit:
\begin{eqnarray}\label{2k1d}
\ad^{2k}fg=[f,\ad^{2k-1}fg]= \\
=(-1)^k\mu_k(q)\left[f, \frac{\partial}{\partial q_k}
\right]=(-1)^{k+1}\mu_k(q)\phi'(q_k-q_{k+1})\frac{\partial}{\partial
p_{k+1}} \ (\mbox{mod }\Lambda^{2k}). \nonumber
\end{eqnarray}

The factor  $\phi'(q_k-q_{k+1})$ at the right-hand side of
(\ref{2k1d}) may vanish at isolated points. Since by induction
hypothesis the vector field $\frac{\partial}{\partial q_k}$ is
tangent to the orbit we can shift our reference point along the
trajectory of $\frac{\partial}{\partial q_k}$ (along the orbit) to a
point where $\phi'(q_k-q_{k+1})$ becomes nonvanishing, while
$\mu_k(q)$ and hence
$-\mu_k(q)\phi'(q_k-q_{k+1})=(-1)^k\prod_{j=1}^{k+1}
\phi'(q_j-q_{j+1})$  remain  nonvanishing. We arrive to a point of
the orbit where the formula (\ref{2k1}) for $\ad^{2k}fg$ and the
formula (\ref{lamk}) for $\Lambda^{2k+1}$ is valid.

The induction step from $\Lambda^{2k+1}$ to $\Lambda^{2k+2}$ can be
accomplished in a similar way. $\Box$

\section{Global controllability of double forced multiparticle system}
\label{dic}

It is easy to see that single-forced multiparticle system  is in
general  uncontrollable, i.e. its attainable sets may not coincide
with the whole state space. For example, if the particle
$\mathcal{P}_n$ is not subject to controlled force, the initial
value $p^0_n$ is positive, and $\phi(q) >0$ (as it is in the case of
Toda lattice), then, given the nature of the interaction between
particles we conclude from the corresponding equation
$\dot{p}_n=\phi(q_{n-1}-q_{n})>0 $, that is $p_n$ is increasing with
time and can not attain values smaller than $p^0_n$.

In this Section and further on we will study the double-input case,
in which controlled forces are applied to the particles
$\mathcal{P}_1$ and $\mathcal{P}_n$.

\subsection{Preliminaries: recurrency of the drift and controllability}
\label{prerecur}

For a control-affine system \eqref{afgen} full dimensionality of its
Lie envelope does not imply in general global controllability. An
obstruction could be actuation of the vector field $f$. It can
provoke a drift in certain direction which can not be compensated by
any control. Now we will formulate conditions under which such
compensation is possible.

 Let the vector field $f$ in $\mathbb{R}^N$ be complete. A point $x
 \in \mathbb{R}^N$ is {\it non-wandering} for $f$ if for
 each its neighborhood $U_x$ and each $t>0$ there exist $x' \in U_x,
 \  t'>t$ such that $e^{tf}(x') \in U_x$. The vector field is
 {\it recurrent} if all the points of $\mathbb{R}^N$ are non-wandering for $f$.

 Theorem due to B.Bonnard and C.Lobry (\cite{Bo},\cite{Lo}) allows to conclude
$\mathcal{A}_{\tilde{x}}=\mathbb{R}^N$ for the system \eqref{afgen}
whenever $\dim \mbox{Lie}_{x}=N$ and the drift vector field is
recurrent.

\begin{prethm}Let $\dim \mbox{Lie}_{x}=N$ and $f$ be recurrent. Then
the system \eqref{afgen} is globally controllable. $\square$
\end{prethm}

\subsection{Global controllability of double forced multiparticle
system}

We consider the same  multi-particle system described by equations
(\ref{dqq}),(\ref{dpq})  but now controlled by forces $u,v$ applied
to the particles $\mathcal{P}_1$ and $\mathcal{P}_n$. The equations
for the momenta of these particles become
\begin{equation}
\dot{p}_1=-\phi(q_{1}-q_{2})+u, \ \dot{p}_n=\phi(q_{n-1}-q_{n})+v,
\label{puv}
\end{equation}
Adjoining these equation to the equations \eqref{dqq}-\eqref{dpq} we
obtain a particular kind of a double-input control-affine  system of
the form
\begin{equation}\label{dicas}
    \dot{x}=f(x)+g^u(x)u+g^v(x)v, \ g^u=\frac{\partial}{\partial
    p_1}, \ g^v=\frac{\partial}{\partial p_n}
\end{equation}
where $f$ is defined by (\ref{vef}).

Our goal is  to prove {\it global controllability} of this system.
To achieve it we  will design the input $u$ as a sum of a certain
smooth feedback control and of an open loop control:
\begin{equation}\label{cu}
u=u_f(q_1)+u_o(t), \ u_f(q_1)=-\frac{\partial U_f}{\partial q_1}.
\end{equation}
We choose the other input $v$ to be a smooth feedback control:
\begin{equation}\label{cv}
v=v_f(q_n)=-\frac{\partial V_f}{\partial q_n}.
\end{equation}

The conditions we impose on $U_f,V_f: \mathbb{R} \rightarrow
\mathbb{R}$ are
\begin{equation}\label{gro}
 \lim_{q \rightarrow -\infty}U_f(q)=+\infty, \
 \lim_{q \rightarrow +\infty}V_f(q)=+\infty
; \end{equation} $U_f,V_f$ are bounded below.

Feeding the controls $\eqref{cu}$ and $\eqref{cv}$ into the
equations  \eqref{puv} we obtain
\begin{equation}
\label{dsis} \dot{p}_1=-\phi(q_{1}-q_{2})+u_f(q_1)+u_o, \
\dot{p}_n=\phi(q_{n-1}-q_{n})+v_f(q_n).
\end{equation}
Now (\ref{dqq})-(\ref{dpq})-(\ref{dsis}) can be treated as a
single-input system with scalar control $u_o$.

Note that we  have proceeded with a particular type of feedback
transformation; more comments on these transformations appear in
Subsection~\ref{sft}.

The drift vector field for the transformed system
\eqref{dqq}-\eqref{dpq}-\eqref{dsis} is Hamiltonian with the
Hamiltonian function
\begin{equation}\label{hamfeed}
   H_{\tilde{f}}(q,p)=\frac{1}{2}\sum_{k=1}^np_k^2+\sum_{j=1}^{n-1}\Phi(q_j-q_{j+1})
   +U_f(q_1)+V_f(q_n).
\end{equation}

Hamiltonian vector fields are recurrent provided that the Lebesgue
sets of the respective  Hamiltonians  functions are compact. Indeed
the Lebesgue sets are invariant for Hamiltonian vector fields, whose
flows are volume-preserving. By Poincare theorem all the
trajectories of such flows must be recurrent.

Therefore it suffices to prove the  following technical lemma.

\begin{lem}
Level sets and Lebesgue sets of the modified Hamiltonian $H_f$ are
compact. $\Box$
\end{lem}

\begin{proof}
Closedness  of the level sets $\{(q,p)| \ H(q,p) = c\}$ and of the
Lebesgue sets $\{(q,p)| \ H(q,p) \leq c\}$ is obvious by the
continuity of $H_f$. It suffices to prove boundedness of the
Lebesgue sets.

Since $\sum_{j=1}^{n-1}\Phi(q_j-q_{j+1})
   +U_f(q_1)+V_f(q_n)$ is bounded below, say by $-B \ (B \geq 0)$,
   then the inequality $H_f \leq c$ implies two constraints:
$$\|p\|^2 \leq c+B \ \bigwedge \
\sum_{j=1}^{n-1}\Phi(q_j-q_{j+1})
   +U_f(q_1)+V_f(q_n) \leq c.$$
Once again by lower boundedness of the functions $\Phi, U_f,V_f$ and
due to  the growth conditions  (\ref{grof}),(\ref{gro}) we conclude
existence of a constant $b$ such that
\begin{eqnarray}
  \sum_{j=1}^{n-1}\Phi(q_j-q_{j+1})
   +U_f(q_1)+V_f(q_n) \leq c \Rightarrow \nonumber \\
  \Rightarrow -q_1 \leq b \bigwedge q_1 -q_2 \leq b
   \bigwedge \cdots q_{n-1}-q_n \leq b \bigwedge q_n \leq b.
\label{qb}
\end{eqnarray}

Summing the first $k$ inequalities in the right-hand side of the
implication (\ref{qb}) we obtain $-q_k \leq kb$, or $q_k \geq -kb$
while summing $n+1-k$ inequalities, starting from the last one, we
obtain $q_k \leq (n+1-k)b$.
\end{proof}

Now we formulate the main result of  this Section.

\begin{thm}
\label{gloc}
 The double-forced multiparticle  system (\ref{dqq})-(\ref{dpq})-(\ref{puv}) is globally
controllable. $\Box$
\end{thm}

\begin{proof}
We invoke controls of the form (\ref{cu})-(\ref{cv}) or,in other
words, aim at  establishing  controllability of the single-input
system (\ref{dqq})-(\ref{dpq})-(\ref{dsis}) controlled by $u_o$.

The Hamiltonian drift vector field $\tilde{f}$
 corresponds to the Hamiltonian (\ref{hamfeed}) with compact Lebesgue
 sets.  By the aforesaid  the drift vector field is
{\it recurrent}, and according to the Proposition~\ref{lier}  the
evaluation of the Lie envelope $\mbox{Lie}_x\{\tilde{f},g^u,g^v\}$
is $2n$-dimensional at every point $x \in \mathbb{R}^{2n}$.

Then by Bonnard-Lobry theorem
 the  single input control-affine system \eqref{dqq}-\eqref{dpq}-\eqref{dsis} is
globally controllable, if the control parameter $u_o$ is allowed to
admit values of both signs: $u_o \in \Omega_o=[-\omega_o, \omega_o],
\ \omega_o>0$.

This implies controllability of the double-input system
(\ref{dqq})-(\ref{dpq})-(\ref{puv}) by means of controls of the form
(\ref{cu})-(\ref{cv}).
\end{proof}

Let us draw conclusions about the constraints, which can be imposed
onto the  values of the controls (\ref{cu}),(\ref{cv}) in order to
keep system controllable. The feedback components $u_f,v_f$ of these
controls are defined via the functions $U_f,V_f$, which can be
chosen globally Lipschitzian with any Lipschitz constant $\omega_o
>0$ in addition to (\ref{gro}). Then the controls
(\ref{cu}),(\ref{cv}) will fit the constraints
$$u_f(t)+u_o(t) \in [-2\omega_o, \omega_o], \ v_f(t) \in  [-\omega_o, \omega_o].$$
It is worth noting that choosing in addition $V_f$ monotonously
increasing we may constrain $v_f$ by the interval $[-\omega_o, 0]$.

We conclude with a Proposition.

\begin{proposition}
\label{prop:3.3} For each $\omega >0$ the two-input system
(\ref{dqq})-(\ref{dpq})-(\ref{puv}) is globally controllable by
means of controls, which meet the constraints
\begin{equation}\label{uvom}
u(t) \in  [-\omega,\omega], \  v(t) \in [-\omega ,0].
\end{equation}
For each pair of points $x^0,x^1$ in the state space of this system,
there exist controls satisfying (\ref{uvom}) which steer the system
from $x^0$ to $x^1$ in some time $T(x^0,x^1,\omega). \ \Box$
\end{proposition}

\begin{remark}
We should mention the publication \cite{PT} where the authors
studied controllability of Toda lattice (in Flaschka form) by means
of $n$-dimensional controls
\begin{eqnarray}
  \dot{a}_1=2b_1^2+u_1, \ \dot{a_2}=2(b_2^2-b_1^2), \ldots , \dot{a}_{n-1}=2(b_{n-1}^2-b_{n-2}^2),
  \dot{a}_n=-2b_{n-1}^2; \nonumber \\
\dot{b}_1=b_1(a_2-a_1)+u_{n+1}, \ldots
\dot{b}_{n-1}=b_{n-1}(a_n-a_{n-1})+u_{2n-1}. \label{Put}
\end{eqnarray}
Note that the controls $u_{n+1}, \ldots u_{2n-1}$  appear in
'kinematic part' of the Toda equations and therefore can not be seen
as forces. There is some controversy (possibly due to  typos) in
what regards  the main result announced in \cite{PT}. The system
\eqref{Put} is not globally controllable on the contrary to what is
claimed, at least because the variable $a_n$ is decreasing according
to \eqref{Put}.
\end{remark}

\section{Feedback linearizability and constant rank}
\label{feelin}

In this Section we demonstrate that double forced multiparticle
system possesses same local properties as controllable linear
system, and in fact is locally equivalent to such a system. To do
this we have to impose an additional regularity assumption onto the
potential of interaction.

\begin{asmp}
\label{asum1} In  Sections \ref{feelin}-\ref{toc} we will  assume
the derivative $\phi'(\cdot)$ of the interaction force $\phi$ to be
nonvanishing.
\end{asmp}

\begin{remark}
This assumption is valid for Toda lattice. $\Box$
\end{remark}

 We start with recalling what are  state-feedback transformations asnd go on with formulation of state-feedback
 linearizability criterion.

\subsection{State-feedback transformation and linearizability}
\setcounter{subsubsection}{1} \label{sft}
 State transformation is a local (at $x^0$)
diffeomorphism $P:x \mapsto y$ of $\mathbb{R}^N$, which acts on the
vector fields of a control-affine system \eqref{afgen} by
differential $P_*$. This results in a state transformation
\begin{equation}\label{st}
\dot{y}=P_*f(y)+\sum_{j=1}^r P_*g^j(y)u_j.
\end{equation}
of \eqref{afgen}.

Feedback transformation is a map
$$v \mapsto u=\alpha(x)+\beta(x)v, \
\beta(x) \mbox{ - nonsingular $(r \times r)$-matrix},$$ where
$\alpha(x),\beta(x)$ are defined in some neighborhood of $x^0$.

 Such a transformation
results in control system $\dot{x}=\bar{f}(x)+\bar{G}(x)v$ with
\begin{equation}\label{fgb}
\bar{f}(x)=f(x)+G(x)\alpha(x), \ \bar{G}(x)=G(x)\beta(x). \ \Box
\end{equation}

\begin{predef}[state-feedback linearizability]
\label{sfl} System is locally state-feedback linearizable if there
exist a local feedback transformation (\ref{fgb}) and a local state
transformation (\ref{st}) such that
\begin{equation}\label{36}
P_* \bar{f}(y)=Ay, \ P_* \bar{G}(y)=B,
\end{equation}
where $A$ is $N \times N$-matrix, $B=\left(b^1 \ldots b^r\right)$
and the vector fields $b^i$ are constant. $\Box$
\end{predef}

\begin{prerem}
On the contrast to standard definition \cite{ASkv,NVS} we do not
require local diffeomorphism $P$ which appears in \eqref{st} and
\eqref{36} maps neighborhood of $x^0$ onto a neighborhood of the
{\it origin}  in $\mathbb{R}^N$. Linearizability means
state-feedback equivalence of the original system to a linear system
defined in a neighborhood of some point $y^0 \in \mathbb{R}^N. \
\Box$
\end{prerem}

We will invoke the following criterion of local state-feedback
linearizabiity which is due to contributions of Jakubczyk-Respondek
and Hunt-Sue-Meyer (\cite{JR},\cite{HSM}).

\begin{prethm}[\cite{JR},\cite{HSM}]
\label{JaR} Smooth control system
$\dot{x}=f(x)+G(x)u=f(x)+\sum_{j=1}^rg^j(x)u_j$ on  $N$-dimensional
state space with $r$-dimensional control $u=(u_1, \ldots , u_r)$ is
locally (at a point $x^0$) state-feedback equivalent to a
controllable linear system, if and only if the vector distributions
\begin{equation}\label{delm}
x \mapsto \Delta^m_x=\spn \{\ad^k f g^j|_{(x)}, k=0, \ldots , m-1; \
j=1, \ldots ,r\}
\end{equation}
possess locally constant dimensions, are involutive, and $\dim
\Delta^n_{x^0}=n. \ \square$
\end{prethm}

\subsection{State-feedback linearizability of double-forced
multiparticle system}

We are going to prove in this Subsection

\begin{thm}
\label{flin}  The double-forced multi-particle system
(\ref{dqq})-(\ref{dpq})-(\ref{puv}) is locally state-feedback
linearizable at each point. $\Box$
\end{thm}

According to the Theorem~\ref{JaR} for establishing state-feedback
linearizability of the double-input control system
(\ref{dqq})-(\ref{dpq})-(\ref{puv}) one has to verify involutivity
of the distributions
\begin{eqnarray}
x \mapsto \Lambda^m_x=\spn \{(\ad^k f g^u)(x), \ k=0, \ldots , m-1\}
\label{deflam} \\
  x \mapsto \Xi^m_x=\spn \{(\ad^k f g^v)(x), \
k=0, \ldots , m-1\}. \label{defxi}\\
\Delta^m=\Lambda^m+\Xi^m.
\end{eqnarray}

Involutivity and constancy of dimensions of these  distributions are
fulfilled by virtue of the following technical lemma.

\begin{lem}
\label{spacoord} For each $k \geq 0$:

i)  the distribution  \eqref{deflam} is constant (does not depend on
$x$); for $m=2k$ and $m=2k+1$
\begin{equation}\label{lamkd}
\Lambda^{2k}=\spn\left\{\frac{\partial}{\partial
p_s},\frac{\partial}{\partial q_s}\left| \  s=1, \ldots ,
k\right.\right\},\
\Lambda^{2k+1}=\Lambda^{2k}+\spn\left\{\frac{\partial}{\partial
p_{k+1}}\right\}
\end{equation}
respectively; $\Lambda^0=\{0\}$.

ii)  the distribution \eqref{defxi}
 is constant (does not depend on $x$); for $m=2k$ and
$m=2k+1$
\begin{equation}\label{delk}
\Xi^{2k}=\spn\left\{\frac{\partial}{\partial
p_s},\frac{\partial}{\partial q_s}\left| \  s=n-k+1, \ldots n\right.
\right\},\ \Xi^{2k+1}=\Xi^{2k}+\spn\left\{\frac{\partial}{\partial
p_{n-k}}\right\}
\end{equation}
respectively; $\Xi^0=\{0\}. \ \Box$
\end{lem}

\begin{corollary}\label{cor42}
The distribution $\Delta^m$ defined by (\ref{delm}) is constant:

Its evaluation at each point coincides with a coordinate subspace
$$q_i=\cdots =q_{i+r}=0, p_j=\cdots =p_{j+s}=0,$$ and obviously is
involutive. Besides $\Delta^{2k}=\mathbb{R}^{2n}$, whenever $2k \geq
n. \ \Box$
\end{corollary}

{\it Proof of the Lemma}. The items i) and ii) are proved in a
similar way;  both proofs follow the course of the proof of
Lemma~\ref{complam}. An additional fact involved is that the factor
$(-1)^k\prod_{j=1}^k \phi'(q_j-q_{j+1})$ which multiplies the vector
field $\frac{\partial}{\partial p_{s+1}}$ in  (\ref{bfp})  is {\it
nonzero} by Assumption~\ref{asum1} at the beginning of the Section,
and therefore (\ref{lamkd}) and (\ref{delk}) are satisfied {\it at
all points}. $\Box$

The conclusion follows  from the Corollary~\ref{cor42}  by
application of the Theo\-rem~\ref{JaR}.

\subsection{Kronecker or controllability indices}

One can draw conclusions on the structure of resulting linear
control double-input system. For linear system {\it Kronecker
indices} form full set of state feedback invariants of a linear
system and determine its Brunovsky normal form.

One can construct sort of Brunovsky normal form for the double
forced multiparticle system; the Kronecker indices are now called
{\it controllability indices} (\cite{NVS}). It turns out that they
depend on whether the number $n$ of particles is even or odd.

For even $n=2\ell$ the two controllability indices  are equal:
$k_1=k_2=n$ , while $k_1= n+1, \ k_2=n-1$ for $n=2\ell-1$. Recall
that the state is $2n$-dimensional.

In both cases we define two sequences of functions by iterated
directional derivation.

For even $n=2\ell$
\begin{eqnarray}\label{flineven}
 y_1=q_{\ell}, \ y_2=L_f y_1 ,  \ldots , \ y_{n}=  L_f y_{n-1}  \\
z_1=q_{\ell+1}, \ z_2=L_f z_1 ,  \ldots , \ z_{n}=  L_f z_{n-1};
\nonumber
\end{eqnarray}
for odd $n=2\ell +1$
\begin{eqnarray}\label{flinodd}
 y_1=q_{\ell+1}, \ y_2=L_f y_1 ,  \ldots , \ y_{n+1}=  L_f y_{n}  \\
z_1=q_{\ell+2}, \ z_2=L_f z_1 ,  \ldots , \ z_{n-1}=  L_f z_{n-2}.
\nonumber
\end{eqnarray}

\begin{lem}
\label{locor}
 The map
$$(q_1, \ldots , q_n,p_1, \ldots , p_n) \mapsto (y_1, \ldots , y_n,
z_1, \ldots , z_n)$$ defined by \eqref{flineven} for $n=2\ell$ and
by \eqref{flinodd} for $n=2\ell+1$ are local diffeomorphisms at each
point of $\mathbb{R}^{2n}$; in both cases $(y_1, \ldots , y_n, z_1,
\ldots , z_n)$ provide a system of local coordinates at each point.
$\Box$
\end{lem}

\begin{thm}
\label{theor43} For  $n$ even,  the double-forced multi-particle
system (\ref{dqq})-(\ref{dpq})-(\ref{puv}) takes  in local
coordinates \eqref{flineven} the form
\begin{eqnarray}\label{39}
  \dot{y}_j=y_{j+1},  \ j=1, \ldots, n-1, \ \dot{y}_n=Y(y,z)+\lambda(y,z)u,;   \\
   \
\ \dot{z}_j=z_{j+1}, \ j=1, \ldots, n-1, \
\dot{z}_n=Z(y,z)+\mu(y,z)v; \ \lambda(y,z)\mu(y,z)\neq 0. \nonumber
\end{eqnarray}
and after a feedback transformation $\bar{u}=Y(y,z)+\lambda(y,z)u, \
\bar{v}=Z(y,z)+\mu(y,z)v,$ the form
\begin{equation}\label{chin}
y_1^{(n)}=\bar{u}, \  z_1^{(n)}=\bar{v}.
\end{equation}

For $n$  odd, the double-forced  multiparticle system takes  in
local coordinates \eqref{flinodd} the form
\begin{eqnarray}\label{40}
  \dot{y}_j=y_{j+1}, \ j=1, \ldots, n;  \ \dot{y}_{n+1}=Y(y,z)+\alpha(y,z)u+\beta(y,z)v,  \\
\dot{z}_j=z_{j+1}, j=1, \ldots, n-2, \
\dot{z}_{n-1}=Z(y,z)+\gamma(y,z)v;\ \alpha(y,z)\gamma(y,z) \neq 0.
\nonumber
\end{eqnarray} and after a feedback transformation
$$\bar{u}=Y(y,z)+\alpha(y,z)u+\beta(y,z)v,\  \bar{v}=\gamma(y,z)v, $$
the form
\begin{equation}\label{chinodd}
y_1^{(n+1)}=\bar{u}, \ z_1^{(n-1)}=\bar{v}. \ \Box
\end{equation}
\end{thm}

\begin{remark}
The construction of the coordinates \eqref{flineven},\eqref{flinodd}
and the linearized forms \eqref{chin}, \eqref{chinodd} of the
controlled multiparticle system are related to  {\it flatness}  of
this system. In particular $y_1,z_1$ cam be seen as flat outputs of
the system. We do not follow this terminology further; addressing
interested readers to the publications \cite{FR,FR2} and references
therein. $\Box$
\end{remark}


\subsubsection{Proofs of Lemma~\ref{locor} and
Theorem~\ref{theor43}} \hspace{3mm}

{\it Proof of Lemma~\ref{locor}} We provide a proof for the odd case
$n=2 \ell +1$. First check that
\begin{eqnarray}\label{adfgu}
\hspace{6mm} L_{\ad^j f g^u} y_r=\left\{
                  \begin{array}{ll}
                    0,  \ j+r < n+1,\\
                    \neq 0,   \ j+r = n+1,
                  \end{array}
                \right.  L_{\ad^j f g^v} y_r=0,  \ j+r <
                n+1.
      \\
\hspace{6mm} L_{\ad^j f g^v} z_r=\left\{
                  \begin{array}{ll}
                    0,   \ j+r < n-1,\\
                    \neq 0,   \ j+r = n-1,
                  \end{array}
                \right.    L_{\ad^j f g^u} z_r=0,  \ j+r \leq n-1.
\label{adfgv}
\end{eqnarray}

We prove the relations \eqref{adfgu} for the coordinates $y_r$ by
induction on $r$. Let $r=1$; then according to the statement i) of
Lemma~\ref{spacoord} $$L_{\ad^j f g^u} y_1=L_{\ad^j f g^u}
q_{\ell+1}=0, \ \mbox{if}\ j < 2\ell +1=n, \ L_{\ad^{n} f g^u }y_1
\neq 0.$$ According to the statement ii) of the same Lemma $L_{\ad^j
f g^v} y_1=0$ for $j < 2\ell +1=n$.

Assuming relations \eqref{adfgu} to be valid for $r < k$, we use the
identity $L_{[f,g]}=L_f \circ L_g-L_f \circ L_g$ to  conclude for
$j+k \leq n+1$
\begin{eqnarray*}
  L_{\ad^j f g^u} y_k= L_{\ad^j f  g^u} L_f  y_{k-1}=-L_{[f,\ad^jf g^u]}y_{k-1} +L_f L_{\ad^jf
  g^u} y_{k-1} =  \\
  = - L_{\ad^jf g^u} y_{k-1}.
  \end{eqnarray*}
We invoked the equality $L_{\ad^jf g^u} y_{k-1}=0$ which is valid by
induction hypothesis. We conclude that $L_{\ad^j f g^u} y_k=-
L_{\ad^jf
  g^u} y_{k-1}$ vanishes, if $j+k < n+1$, and is different from $0$ if $j+k =
  n+1$. Similar reasoning settles the induction passage for $L_{\ad^j f
  g^v} y_k$. We proceed along the same lines in the proof of
  \eqref{adfgv}.

The differentials $dy_1, \ldots , dy_{n+1}, dz_1, \ldots , dz_{n-1}$
are dual to {\it linear independent system} of the vector fields
$$\ad^jfg^u, \ j=0, \ldots, n; \ \ad^i fg^v, \ i=0, \ldots, n-2.  $$
Hence \eqref{flinodd} defines local coordinate system in the odd
case. Proof for the even case is similar. $\Box$

{\it Proof of Theorem~\ref{theor43}}. Again the proofs of \eqref{39}
and \eqref{40} are similar; we sketch the second one.

First according to \eqref{adfgu}
$$L_{g^u}y_r=0,\ L_{g^v}y_r=0, \ \mbox{for} \  r \leq
n, \ L_{g^v }z_s=0, \ \mbox{for} \ s <n-1, \ L_{g^u} z_s=0, \
\mbox{for} \ s \leq n-1.$$  Also $L_{g^u} y_{n+1} \neq 0, \ L_{g^v}
z_{n-1} \neq 0$.

Basing on these identities we compute $$\dot{y}_j=L_{(f + g^u u+g^v
v)} y_j=L_f y_j=y_{j+1}$$ for $j < n+1$. Also $$\dot{z}_s=L_{(f +
g^u u+g^v v)}z_s=L_f z_s=z_{s+1}$$ for $s < n-1$.

By the same computation
\begin{eqnarray*}
  \dot{y}_{n+1}=L_{(f + g^u u+g^v v)}
y_{n+1}=L_f y_{n+1}+ (L_{g^u} y_{n+1})u + (L_{g^v} y_{n+1})v. \\
 \dot{z}_{n-1}=L_{f + g^u u+g^v v}
z_{n-1}=L_f z_{n-1}+ (L_{g^v} z_{n-1})v,
\end{eqnarray*}
where $L_{g^u} y_{n+1},L_{g^v} z_{n-1}$ are nonvanishing functions
of $y_i,z_j. \ \Box$

\subsection{Systems of constant rank}
\label{conrank}

 \setcounter{subsubsection}{1}

  We will discuss another property of
the control system (\ref{dqq})-(\ref{dpq})-(\ref{puv}) which follows
from its state-feedback linearizability. It is called {\it constancy
of rank} and has been introduced by A.Agrachev and S.A.Vakhrameev in
\cite{AVah}.

\begin{predef}
For a control system $\dot{x}=f(x,u)$ consider input/end-point map
$\mathcal{E}_{x^0,T}$ (with $x^0, T>0$ being  parameters) which puts
into correspondence to each admissible control (input) $u(\cdot)$
the point $x(T)$ of the corresponding trajectory of the control
system starting at $x^0$ at $t=0$ and driven by the control
$u(\cdot)$. We denote the map $\mathcal{E}_{x^0,T}(u(\cdot))$. The
system $\dot{x}=f(x,u)$ is of constant rank if for each $x^0,T$ the
rank (the differential) of $\mathcal{E}_{x^0,T}(u(\cdot))$ does not
depend on $u(\cdot).\ \Box$
\end{predef}

The systems of constant rank inherit many properties of linear
systems. It is known that state-feedback linearizable systems
possess constant rank.

\begin{corollary}
The controlled double forced multiparticle system
(\ref{dqq})-(\ref{dpq})-(\ref{puv}) possesses constant rank. $\Box$
\end{corollary}

\section{Time-optimal control for double-forced multi-particle system}
\label{toc}

Let us consider a problem of time-optimal relocation of  particles
of the double-forced multiparticle system described by  equations
(\ref{dqq})-(\ref{dpq})-(\ref{puv}) with control parameters
constrained by (\ref{uvom}).

\begin{probl}
\label{oprel} Given two points
$\tilde{x}=(\tilde{q},\tilde{p}),\hat{x}=(\hat{q},\hat{p})$ (two
couples of initial and final values of positions and momenta of the
particles) find a pair of admissible controls which  steer the
system (\ref{dqq})-(\ref{dpq})-(\ref{puv})-\eqref{uvom} from
$\tilde{x}$ to $\hat{x}$ in a minimal time $T>0$. $\ \Box$
\end{probl}

Existence of an optimal control in a control-affine problem with
bounded convex set of control parameters follows from  Filippov's
Theorem (\cite{Ces}).

We will be interested in structure of optimal controls and start
with formulation of Pontryagin maximum principle - necessary
optimality condition for time optimal control problem. We limit
ourselves to control-affine problems.

\subsection{Time-optimal control, Pontryagin Maximum Principle, bang-bang extremals}
\setcounter{subsubsection}{1}

Consider time-optimal control problem under boundary conditions
\begin{equation}\label{bocon}
x(0)=\tilde{x}, \ x(T)=\hat{x}, \ T \rightarrow \min .
\end{equation}
 for control-affine system
\eqref{afgen}. One seeks an admissible control, which steers the
system \eqref{afgen} from $\tilde{x}$ to $\hat{x}$ in minimal time
$T$. Along this Subsection the set $U$ of control parameters in
\eqref{afgen} is assumed
  to be a compact convex polyhedron in
$\mathbb{R}^r$.

A first-order necessary condition for $L_1$-local optimality of an
admissible control $\tilde{u}(\cdot)$ for such problem is provided
by {\em Pontryagin Maximum Principle}
 (see \cite{Pon}).

\begin{prethm}
Let  pair $(\tx ,\tu)$ be a minimizing control and corresponding
trajectory for the time-optimal control problem
\eqref{afgen}-\eqref{bocon}, $T$ being minimal time. Then there
exists a non-zero absolutely continuous covector-function $\tps : R
\rightarrow \left(\mathbb{R}^N\right)^*$, such that the pair
$(\tx,\tps)$ satisfies Hamiltonian system with the Hamiltonian:
\begin{equation}                                    \label{ham}
\Pi(x,\psi ,u)= \langle \psi ,f(x)\rangle+ \langle \psi , G(x)u
\rangle.
\end{equation}
In local coordinates this system takes form
$$\dot{x}=\partial \Pi/ \partial \psi (\txe, \psi , \tue(\tau))
                   , \  \dot{\psi}=-\partial \Pi/
\partial x (\txe, \psi , \tue(\tau)).$$
Besides the following conditions hold:

i) Maximality Condition:
\begin{equation}                      \label{maco}
\Pi(\txe(t), \tps(t), \tue(t))= \max \{\Pi(\txe(t), \tps(t),u): \ u
\in U \} \ \mbox{a.e. on} \ [0,T];
\end{equation}

ii) Transversality Condition:
$$\Pi(\txe(T),\tps(T),\tue(T)) \geq 0. \ \Box$$
\end{prethm}


The solutions of the equations of the Pontryagin Maximum
 Principle are called {\em Pontryagin extremals}, the corresponding
 controls $\tu$ are called {\em extremal controls}.

For any $\tau \in [0,T]$ maximum \eqref{maco} of the control-affine
Hamiltonian \eqref{ham}, is attained at some  face of the polyhedron
$U$. This face can be $0-$dimensional,  then extremal control takes
its value at a vertex of the polyhedron,  or $s$-dimensional ($0<s
\leq r$) and then the maximality condition does not determine the
value of extremal control uniquely.

We  call {\it bang-bang} the extremal controls for which the maximum
is achieved  at some  vertices of the polyhedron $U$ on a set of
full measure in $[0,T^*]$. Change of the value of control from one
vertex to another one is called {\it switching}. The controls which
take their values on faces of positive dimensions are called {\it
singular}; it will turn out that such controls do not occur in our
problem .

A classical bang-bang result for linear time-optimal control problem
with the dynamics $\dot{x}=Ax+Bu, \ u \in U$ the following theorem
on structure of optimal controls  has been proven by R.V.Gamkrelidze
(see \cite{Pon}).

\begin{preprop}
\label{Gamba} If for a directing vector $V$ of any edge of the
polyhedron $U$ the vectors $BV,ABV, \ldots , A^{n-1}BV$ are linearly
independent (genericity assumption), then any control which
satisfies the Pontryagin maximum principle (and in particular any
optimal control) is piecewise constant, takes its values at the
vertices of the polyhedron $U$ and possesses finite number of
switchings. $\ \Box$
\end{preprop}

{\it Nonlinear} control-affine time-optimal problem
\eqref{afgen}-\eqref{bocon} do not resemble in general linear
time-optimal problems and in particular the conclusion of the
Proposition~\ref{Gamba} does not hold  for them generically.

In the next  Subsection we will prove that controls providing
time-optimal relocation of particles (Problem 1) are bang-bang.

\subsection{Time-optimal relocation problem. Bang-bang properties of optimal controls}
\setcounter{subsubsection}{1}

The key additional feature of the control system
\eqref{dqq}-\eqref{dpq}-\eqref{puv} which allows to establish the
bang-bang property is its constancy of rank and feedback
linearizability (see Section~\ref{feelin}).

\begin{thm}
\label{baba} Optimal controls for the time-optimal relocation
problem  (Problem 1) are bang-bang and possess finite number of
switchings. $\Box$
\end{thm}

The proof is based on a criterion due to A.Agrachev and S.Vakhrameev
(\cite{AVah},\cite{Vah}). The criterion is formulated for the
control-affine time-optimal problem \eqref{afgen}-\eqref{bocon} and
involves the following assumptions.

{\it Genericity Assumption.} For a directing vector $w$ of each edge
of the polyhedron $U$ and for all $x \in \mathbb{R}^N$ the vectors
\begin{equation}\label{genas}
Gw|_x, \ad \ f Gw|_x, \ldots ,  (\ad \ f)^{N-1} Gw|_x,
\end{equation}
 are
linearly independent.

{\it  Bang-bang condition}  is satisfied for an edge $w$ of the
polyhedron $U$ if for each point $\hat{x} \in \mathbb{R}^N$ there
exist smooth covector-functions $x \mapsto a^i_j(x) \in
\mathbb{R}^{r^*}$ defined in some neighborhood $\Omega$ of $\hat{x}$
such that for any $u \in U$ and for all $i=0,1, \ldots $
\begin{equation}\label{gugw}
\left.\left[Gu, (\ad \ f)^i Gw \right]\right|_x= \sum_{j=1}^i
\left.\langle a^i_j(x), u \rangle (\ad \ f)^j Gw\right|_x.\ \Box
\end{equation}

\begin{prethm}[\cite{Vah}]
\label{avah} Let \eqref{afgen} be  analytic \footnote{Actually less
restrictive condition of finite-definiteness is needed.} system of
constant rank, which satisfies the genericity assumption and the
bang-bang condition for each edge of the polyhedron $U$ of
admissible control parameters. Then any time-optimal control of the
problem \eqref{afgen}-\eqref{bocon} is bang-bang with a finite
number of switchings. $\Box$
\end{prethm}

To apply the criterion provided by this Theorem to Problem 1  we
first note that the dynamics of multi-particle system is analytic.
The control system  \eqref{dqq}-\eqref{dpq}-\eqref{puv} is locally
state-feedback linearizable  and hence is of constant rank (equal to
$2n$).

The polyhedron $U$ defined by (\ref{uvom}) is a rectangle. The
directing vectors $w_u,w_v$ of its edges are parallel to the axes
$u$ and $v$. Substituting these vectors in place of $w$ in
(\ref{genas}) we obtain two sequences of vector fields $$g^u|_x,
(\ad \ f) g^u|_x, \ldots ,  (\ad \ f)^{N-1} g^u|_x; \ \mbox{and} \
g^v|_x, (\ad \ f) g^v|_x, \ldots ,  (\ad \ f)^{N-1} g^v|_x,$$ both
of which are linearly independent according to Lemma~\ref{spacoord}.

The validity of the bang-bang condition (\ref{gugw}), verified for
the rectangle $U$, can be derived  from  the equalities
\begin{equation}\label{grogsi}
\left.\left[g^\rho, (\ad \ f)^i g^\sigma \right]\right|_x=
\sum_{j=1}^i \left. a^{i\rho}_{j\sigma}(x) (\ad \ f)^j
g^\sigma\right|_x ,
\end{equation}
where the symbols $\rho$ and $\sigma$ coincide with either $u$ or
$v$.

All these equalities can be verified in a similar way. We do it for
$\rho=v, \sigma=u$, and distinguish the cases of even and odd $i$.

According to the Lemma~\ref{spacoord} , one gets for  $i=2k-1$:
 $$ (\ad \ f)^{2k-1} g^u=\sum_{s=1}^k\alpha_s(x)\frac{\partial}{\partial q_s}+
\beta_s(x)\frac{\partial}{\partial p_s}.$$ As far as
$g^v=\frac{\partial}{\partial p_n}$ is constant vector field, which
commutes with $\frac{\partial}{\partial q_s},
\frac{\partial}{\partial p_s}$, then
$$\left[g^v,(\ad \ f)^{2k-1} g^u\right]=\sum_{s=1}^k\left(L_{g^v}\alpha_s(x)\right)
\frac{\partial}{\partial q_s}+
\left(L_{g^v}\beta_s(x)\right)\frac{\partial}{\partial p_s}.$$ The
values of this Lie bracket belong to  $\Lambda^{2k}$ defined by
(\ref{lamkd}). According to Lemma~\ref{spacoord} this Lie bracket
can be represented as a linear combination $\sum_{j=1}^i
\left.b^i_j(x)(\ad \ f)^j g^u\right|_x.$

The proof for $i=2k$ is obtained similarly.

\subsection{Uniform boundedness of the number of switchings}
\setcounter{subsubsection}{1}

 It is known that for a bang-bang
control the number of switchings can be arbitrarily large and even
infinite.  In this subsection we wish to establish a stronger
property of bang-bang optimal controls for time-optimal relocation
problem.  It guarantees uniform boundedness of the number of
switchings for all optimal trajectories contained in some compact of
$\mathbb{R}^{2n}$.

For control-affine system {\it with single input} the problem has
been formulated  and studied by A.J.Krener (\cite{Kre}) and
H.J.Sussmann (\cite {Sus}).

\begin{definition} Control problem possesses strong bang-bang
property with bounds on the number of switchings, if for every
compact set $K$ and $T>0$ there exists an integer $N(K,T)$ such that
 any time-optimal trajectory of time duration $T$, which connects
two points $\tilde{x},\hat{x}$ and is contained in $K$ is  bang-bang
trajectory with at most $N$ switchings.
\end{definition}

From the aforesaid we already know that all optimal controls in
time-optimal relocation problem are bang-bang. Technical lemma
proved in \cite{Sus} for single-input case can be adapted to the
time-optimal problem (\ref{dqq})-(\ref{dpq})-(\ref{puv}) {\it with
two inputs} given the special Lie structure of the controlled
multiparticle system and the fact that the set (\ref{uvom}) of
control parameters is a rectangle.

\begin{thm}
\label{bbns} Time optimal particle relocation problem for double
forced multiparticle system possesses strong bang-bang property with
bound on the number of switchings. $\Box$
\end{thm}

\begin{proof}

For the  control system  (\ref{dicas})  the Hamiltonian of the
Pontryagin maximum principle takes form
\begin{eqnarray}\label{hapon}
    \Pi(q,p,\psi_q,\psi_p,u,v)=\sum_{k=1}^n\psi_{q_k}p_k+\sum_{\ell=1}^n\psi_{p_\ell}
    \left(\phi(q_{\ell-1}-q_\ell)-\phi(q_{\ell}-q_{\ell+1})\right)+ \nonumber \\
    +\psi_{p_1}u+\psi_{p_n}v.
\end{eqnarray}

 According to (\ref{hapon}) the bang-bang  values of the
controls calculated from the maximality condition (\ref{maco}) are
defined by the sign of the "switching functions"
$\sigma^u(t)=\psi_{p_1}(t) \ \sigma^v(t)=\psi_{p_n}(t)$:
$$u(t)=\omega_0\mbox{sign}\sigma^u(t), \
v(t)=\frac{\omega_0}{2}\left(1+\mbox{sign}\sigma^u(t)\right).$$

Then it suffices to prove that extremal trajectories contained in
any fixed compact $K$ the number of zeros of the switching functions
$\sigma^u(t), \ \sigma^v(t)$ is bounded by a constant $C(K)$.

To this end we  introduce the functions
\begin{equation}\label{siro}
\sigma^\rho_k(t)=\langle \psi(t), (\ad f)^kg^\rho \rangle, \
\sigma^\rho_0(t)=\sigma^\rho(t), \  \rho \in \{u,v\}.
\end{equation}

Evidently
\begin{eqnarray}
\label{sigad}
\dot{\sigma}^\rho_k(t)= \langle \psi(t), (\ad
f)^{k+1}g^\rho(x(t))\rangle+ \\
+u\langle \psi(t), [g^u,(\ad f)^{k}g^\rho](x(t))\rangle+v\langle
\psi(t), [g^v,(\ad f)^{k}g^\rho](x(t))\rangle, \nonumber
\end{eqnarray}
 and by (\ref{grogsi})
\begin{eqnarray*}
  \dot{\sigma}^\rho_k(t)=\sigma^\rho_{k+1}(t)+\sum_{j=1}^k u(t)\alpha^{ku}_{j\rho}(t)\sigma^\rho_{j}(t)+
\sum_{j=1}^k v(t)
\alpha^{kv}_{j\rho}(t)\sigma^\rho_{j}(t)= \\
 = \sum_{j=1}^k
  a_{kj}(t)\sigma^\rho_j(t)+\sigma^\rho_{k+1}(t).
\end{eqnarray*}

 As far as $\{(\ad f)^kg^u | \ k=0, \ldots , 2n-1\}$ span $R^{2n}$,
the iterated Lie bracket  $(\ad f)^{2n}g^\rho$ can be represented as
\begin{equation}\label{2nprev}
(\ad f)^{2n}g^\rho=\sum_{j=0}^{2n-1}\gamma_j(x)(\ad f)^{j}g^\rho.
\end{equation}
Setting  $k=2n-1$ in \eqref{sigad} and substituting \eqref{2nprev}
into its right-hand side we conclude
\begin{equation}\label{last}
\dot{\sigma}^\rho_{2n-1}=\sum
_{j=0}^{2n-1}a_{2n-1,j}(x)\sigma^\rho_{j}.
\end{equation}

Hence the functions $\sigma^\rho_k(t), \ k=0, \ldots , 2n-1$ satisfy
the following quasitriangular system of linear differential
equations
\begin{equation}\label{trian}
  \dot{\sigma}^\rho_k(t)=\sum_{j=1}^k
  a_{kj}(t)\sigma^\rho_j(t)+\sigma^\rho_{k+1}(t), \ k=0, \ldots ,
  2n-2,
\end{equation}
completed by the equation (\ref{last}).

We now apply to the quasitriangular system the following technical
result due to H.J.Sussmann (\cite{Sus}); it establishes a bound on
the number of zeros of the switching function.

\begin{preprop}
\label{telem}  Let absolute values of all the coefficients at the
right-hand side of (\ref{trian})-(\ref{last}) be bounded by a
constant $A>0$. Then there exists positive $T(A)$ such that on any
time interval $\mathcal{I}$ of length $\leq T$ the component
$\sigma^\rho_0(t)$ of the solution of (\ref{trian}) either vanishes
identically, or possesses at most $2n-1$ zeros. $\Box$
\end{preprop}

The component $\sigma^\rho_0(t)$ can not vanish identically on an
interval, as long as then all $\sigma^\rho_k(t)$ must vanish by
virtue of (\ref{trian}), which in its turn is impossible due to the
definition (formula (\ref{siro})) of $\sigma^\rho_k(t)$ and
Proposition~\ref{lier} by which at each point $\dim \spn\{(\ad
f)^kg^\rho| \ k=0, \ldots , 2n-1\}=2n.$

The conclusion of the Theorem~\ref{bbns} follows now from
Lemma~\ref{telem} by a standard reasoning provided in \cite{Sus}; an
additional component needed for the proof is that the bound $A$ for
the coefficients (\ref{trian})-(\ref{last}) can be chosen the same
for all extremal trajectories contained in a compact $K \subset
\mathbb{R}^{2n}.$

\end{proof}


\end{document}